\documentclass[12pt,twoside, epsf]{article}

\usepackage{color,graphicx,times}

\usepackage{amssymb}
\usepackage{amsmath}
\usepackage{amscd}
\usepackage{float}
\usepackage{graphicx}
\usepackage{amsfonts}
\usepackage{pb-diagram}
\usepackage{eufrak}

 \newcommand{\CGlyph}{\raisebox{-0.25\height}{\includegraphics[width=0.5cm]{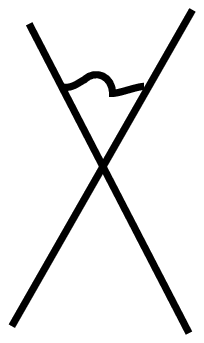}}}
\newcommand{\VGlyph}{\raisebox{-0.25\height}{\includegraphics[width=0.5cm]{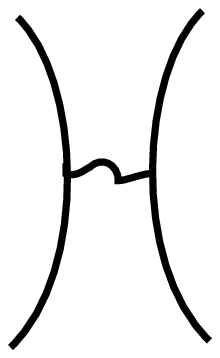}}}
\newcommand{\YGlyph}{\raisebox{-0.25\height}{\includegraphics[width=0.5cm]{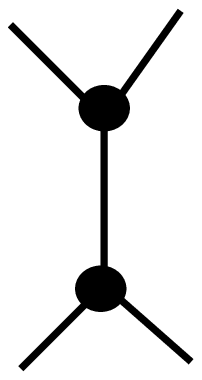}}}

\newcommand{\CDiag}{\raisebox{-0.25\height}{\includegraphics[width=0.5cm]{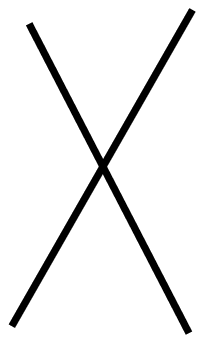}}}
\newcommand{\VDiag}{\raisebox{-0.25\height}{\includegraphics[width=0.5cm]{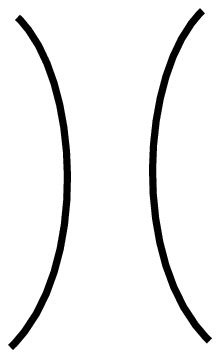}}}
\newcommand{\CDotDiag}{\raisebox{-0.25\height}{\includegraphics[width=0.5cm]{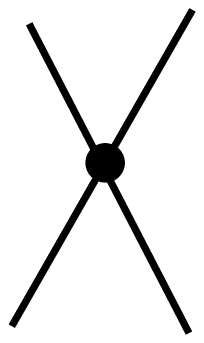}}}
\newcommand{\CCircleDiag}{\raisebox{-0.25\height}{\includegraphics[width=0.5cm]{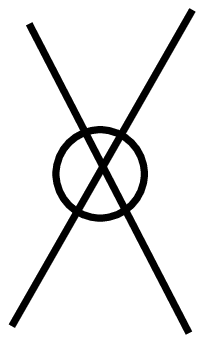}}}

\makeatother

\let \ttorg \tt \def \tt{\ttorg \obeyspaces}

\begin{document}

\date{}

\title{\bf A State Calculus for Graph Coloring}

\author{Louis H. Kauffman \\ 
Department of Mathematics, Statistics and Computer Science \\ 
University of Illinois at Chicago \\ 
851 South Morgan Street\\
Chicago, IL, 60607-7045}

\maketitle

\thispagestyle{empty}

\begin{abstract}
This paper discusses reformulations of the problem of coloring plane maps with
 four colors.  We give a number of alternate ways to formulate the coloring problem including a tautological expansion similar to the Penrose Bracket, and 
 we give a simple extension of the Penrose Bracket that counts colorings of arbitrary cubic graphs presented as immersions in the plane.
 \end{abstract}

\section{Introduction}
This work involves a rewriting of the
coloring problem in terms of two-colored systems of Jordan curves in the plane. 
These systems, called {\em formations} \cite{SB}, are in one-to-one correspondence with
cubic plane graphs that are colored with three edge colors so that three distinct
colors are incident to each node of the graph. It has long been known that the
four color problem can be reformulated in terms of coloring such cubic graphs.
The present paper is an extension of \cite{VCP,LK1,LK2,LK3,LK4} but is self-contained.\\

This paper is organized into five sections:
\begin{enumerate}
\item  Introduction.
\item Cubic Graphs and Formations.
\item Cubic Graphs and the Four Color Problem.
\item Loops and States.
\item The Penrose Formula and a Generalization to All Cubic Graphs.
\end{enumerate}

The second section on cubic graphs and formations shows how we can reformulate the coloring problem in terms of interactions of Jordan curves in the plane. This section is based on previous papers
 \cite{LK1,LK2,LK3,LK4} of the author. We use the fact that cubic graphs have perfect matchings 
to accomplish this aim. In sections  3 and 4 we show how to create tautological state sums for the colorings of arbitrary cubic maps. In working with these states we find a simple example of a plane graph where only one state admits the coloring. This example, being close to being uncolorable, is a good example to modify. We find a natural non-planar modification of it that yields the Isaacs graph $J_{3}$ and hence provides a very neat motivation for the $J$-constructions of Rufus Isaacs \cite{Isaacs}. We give other examples showing how the tautological state sum construction is related to the Petersen graph and other aspects of uncolorability. In Section 5 we prove, using formations, the properties of the well-known Penrose Bracket that counts the number of proper edge colorings of a cubic plane graph. We give a very simple extension of the method for the Penrose Bracket that allows it to compute the number of colorings for an arbitrary cubic graph, presented as an immersed graph in the plane.\\ 
 
\noindent {\bf Acknowledgement.}    It gives the author pleasure to thank
James Flagg for helpful conversations in the course of 
constructing this paper.
\bigbreak

\section{\bf Cubic Graphs and Formations}

A {\em graph} consists in a node set $V$ and an edge set $E$ such that every
edge has two nodes associated with it (they may be identical).  If a node is
in the set of nodes associated with an edge, we say that this node {\em
belongs} to that edge. If two nodes form the node set for a given edge we
say that that edge {\em connects} the two nodes (again the two may be
identical).  A {\em loop} in a graph is an edge whose node set has cardinality
one. In a {\em multi-graph} it is allowed that there may be a multiplicity of
edges connecting a given pair of nodes. All graphs in this paper are
multi-graphs, and we shall therefore not use the prefix ``multi" from here on.\\

A {\em cubic graph} is a graph in which every node either belongs to three
distinct edges, or there are two edges at the node with one of them a loop. A
{\em coloring} (proper coloring) of a cubic graph $G$ is an assignment of the
labels $r$ (red), $b$ (blue), and $p$ (purple) to the edges of the graph so that
three distinct labels occur at every node of the graph. This means that there
are three distinct edges belonging to each node and that it is possible to
label the graph so that three distinct colors occur at each node. Note that a
graph with a loop is not colorable. \vspace{3mm}

The simplest uncolorable cubic graph is illustrated in Figure~\ref{peter}.  For obvious
reasons, we refer to this graph as the {\em dumbell}.  Note that the dumbell is
planar. Figure~\ref{peter} also illustrates a more complex dumbell and the Petersen graph, a non-planar uncolorable.\\

An edge in a connected plane graph is said to be an {\em isthmus} if the deletion of that
edge results in a disconnected graph. It is easy to see that a connected plane
cubic graph without isthmus is loop-free. \vspace{3mm}

Heawood (see \cite{Kempe,Heawood,Petersen}) reformulated the four-color conjecture (which we will henceforth refer to
as the {\em Map Theorem}) for plane maps to a corresponding statement about the
colorability of plane cubic graphs. In this form the theorem reads \\

\noindent {\bf Map Theorem for Cubic Graphs.}  A connected plane cubic graph without
isthmus is properly edge-colorable with three colors. \vspace{3mm}

We now introduce a diagrammatic representation for the coloring of a cubic graph.
Let $G$ be a cubic graph and let $C(G)$ be a coloring of $G.$ Using the colors
$r$, $b$ and $p$ we will write purple as a formal product of red and blue: $$p =
rb.$$

\noindent One can follow single colored paths on the coloring $C(G)$ in the
colors red and blue. Each red or blue path will eventually return to its starting
point, creating a circuit in that color. The red circuits are disjoint from one
another, and the blue circuits are disjoint from one another. Red circuits
and blue circuits may meet along edges in $G$ that are colored purple
($p=rb$).   In the case of a plane graph $G$, a meeting of two circuits may take
the form of one circuit crossing the other in the plane, or one circuit may share
an edge with another circuit, and then leave on the same side of that other circuit.
We call these two
planar configurations a {\em cross} and a {\em bounce} respectively. \vspace{3mm}

\begin{figure}[htb]
     \begin{center}
     \begin{tabular}{c}
     \includegraphics[width=6cm]{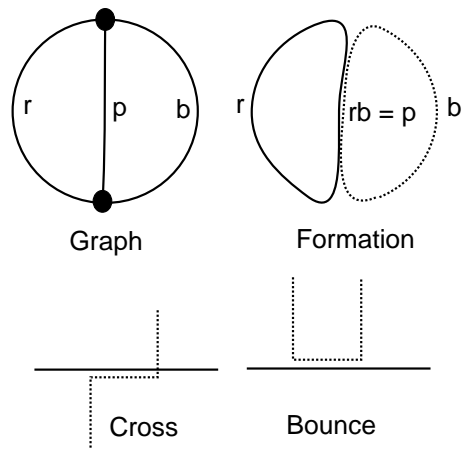}
     \end{tabular}
     \caption{\bf Coloring and Formation.}
     \label{theta}
\end{center}
\end{figure}

\noindent {\bf Definition.} A {\em formation} \cite{SB} is a finite collection of 
simple closed curves in the plane, with each curve colored either red or blue such that the
red curves are disjoint from one another, the blue curves are disjoint from one
another and red and blue curves can meet in a finite number of segments (as
described above for the circuits in a coloring of a cubic graph). \vspace{2mm}

\noindent Associated with any formation $F$ there is a well-defined cubic graph
$G(F)$,  obtained by identifying the shared segments in the formation as edges in
the graph, and the endpoints of these segments as nodes. The remaining
(unshared) segments of each simple closed curve constitute the remaining edges of
$G(F).$  A formation $F$ is said to be a formation for a cubic graph $G$ if $G =
G(F).$ We also say that $F$ {\em formates} $G.$ \vspace{2mm}

\noindent A {\em plane formation} is a formation such that each simple closed
curve in the formation is a Jordan curve in the plane. For a plane formation,
each shared segment between two curves of different colors is either a bounce or
a crossing (see above), that condition being determined by the embedding of the
formation in the plane. \vspace{3mm}

Since the notion of a formation is abstracted from the circuit decomposition of a
colored cubic graph, we have the proposition: \vspace{3mm}

\noindent {\bf Proposition.}  Let $G$ be a cubic graph and $Col(G)$ be the set of
colorings of $G$.  Then $Col(G)$ is in one-to-one correspondence with the set of
formations for $G.$ 
\bigbreak

\noindent In particular, the Map Theorem is equivalent to the
\bigbreak

\noindent {\bf Formation Theorem.}  Every connected plane cubic graph without
isthmus has a formation. \vspace{3mm}

\begin{figure}[htb]
     \begin{center}
     \begin{tabular}{c}
     \includegraphics[width=6cm]{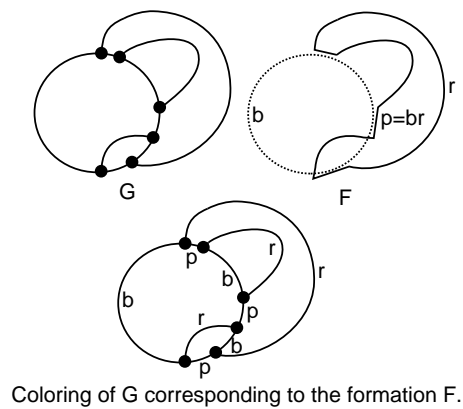}
     \end{tabular}
     \caption{\bf Second Example of Coloring and Formation.}
     \label{colform}
\end{center}
\end{figure}

This equivalent version of the Map Theorem is due to G. Spencer-Brown \cite{SB}. 
The advantage of the Formation Theorem is that, just as one can enumerate graphs,
one can enumerate formations.  In particular, plane formations are generated by
drawing systems of Jordan curves in the plane that share segments according to
the rules explained above. This gives a new way to view the evidence for the Map
Theorem, since one can enumerate formations and observe that all the plane cubic
graphs are occurring in the course of the enumeration!   See Figures~\ref{theta} and \ref{colform} for
illustrations of the relationship of formation with coloring. \vspace{3mm}

\noindent {\bf Remark.}    In depicting formations, we have endeavored to keep the
shared segments slightly separated for clarity in the diagram. These separated
segments are amalgamated in the graph that corresponds to the formation.
\vspace{3mm}

\section{Cubic Graphs and the Four Color Problem}
In this section we state results equivalent to the Map Theorem in terms of perfect matchings in the graph.
We shall use these concepts in conjunction with the notion of formation from the previous section in the rest of the paper.\\

Recall that a cubic graph $G$ is said to be {\em properly colored} with $3$ colors if the edges of $G$
are colored from the $3$ colors so that all colors incident to any node of $G$ are distinct, and that the Map Theorem states
that an isthmus-free planar cubic graph can be properly colored.\\

We shall first use this result to give yet another (well-known) equivalent version of the Map Theorem. 
To this end, call a disjoint collection $E$ of edges of $G$ that includes all the nodes of $G$ a {\em perfect matching}
of $G.$ Then $C(E,G) = G - Interior(E)$ is a collection of {\em cycles} (graphs homeomorphic to the circle, with two
edges incident to each node). We say that $E$ is an {\em even} perfect matching of $G$ if every cycle in $C(E)$ has
an  even number of edges.
\bigbreak

\noindent {\bf Theorem.} The following statement is equivalent to the Map Theorem: Let $G$ be a plane cubic graph with
no isthmus. There there exists an even  perfect matching of $G.$
\bigbreak

\noindent {\bf Proof.} 
Let $G$ be a cubic plane graph with no isthmus. Suppose that $G$ is properly $3$-colored from the set
$\{ a,b,c \}$. Let $E$ denote all edges in $G$ that receive the color $c$. Then, by the definition of proper coloring, the 
edges in $E$ are disjoint. By the definition of proper $3$-coloring every node of $G$ is in some edge of $E.$
Thus $E$ is a perfect matching of $G.$ Since each cycle in $C(E,G)$ is two-colored by the the set $\{ a,b \},$
each cycle is even. Hence $E$ is an even perfect matching of $G.$\\

Conversely, suppose that $E$ is an even perfect matching of $G.$ Then we may assign the color $c$ to all the edges of
$E,$ and color the cycles in $C(E)$ using $a$ and $b$ (since each cycle is even). The result is a proper
$3$-coloring of  the graph $G.$ This completes the proof of the Theorem. $\hfill\Box$ \\

\begin{figure}[htb]
     \begin{center}
     \begin{tabular}{c}
     \includegraphics[width=6cm]{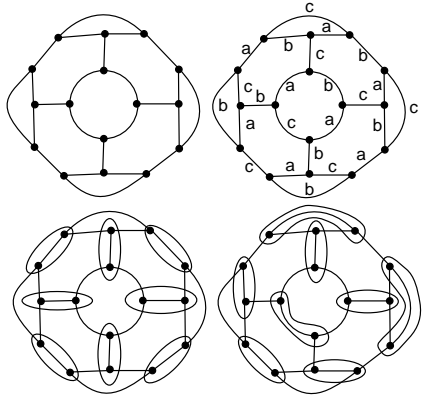}
     \end{tabular}
     \caption{\bf Perfect Matchings of a Cubic Plane Graph}
     \label{perfect}
\end{center}
\end{figure}

\noindent {\bf Remark.} See Figure~\ref{perfect} for an illustration of two perfect matchings of a graph $G.$ One perfect matching is not
even. The other perfect matching is even, and the  corresponding coloring is shown. This Theorem shows that one could
conceivably divide the proving of the FCT into two steps: First prove that every cubic plane isthmus-free graph has a perfect
matching. Then prove that it has an  even  perfect matching. In fact, the existence of a perfect matching is hard, but available
\cite{Petersen}, while the existence of an even perfect matching is really hard!\\

\noindent {\bf Proposition.} Every cubic graph with no isthmus has a  perfect matching.
\bigbreak

\noindent {\bf Proof.}  See  \cite{Petersen}, Chapter $4.$ $\hfill\Box$
\bigbreak

\begin{figure}[htb]
     \begin{center}
     \begin{tabular}{c}
     \includegraphics[width=6cm]{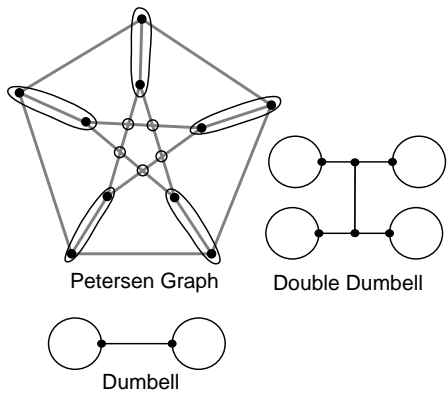}
     \end{tabular}
     \caption{\bf Petersen and Dumbells}
     \label{peter}
\end{center}
\end{figure}

\noindent {\bf Remark.}  There are graphs that are
uncolorable. Two famous such culprits are indicated in Figure~\ref{peter}. These are examples of graphs with  perfect matchings, but
no even perfect matching. The second example in  Figure~\ref{peter} is the ``dumbell graph". It is planar, but has an isthmus. The first
example is the Petersen Graph. This graph is non-planar. We have illustrated the Petersen with one  perfect matching that has two 5-cycles. No  perfect matching of the Petersen is even. The third  "double dumbell" graph illustrated in Figure~\ref{peter} has no
perfect matching.
\bigbreak

\section{Loops and States}
\begin{figure}[htb]
     \begin{center}
     \begin{tabular}{c}
     \includegraphics[width=6cm]{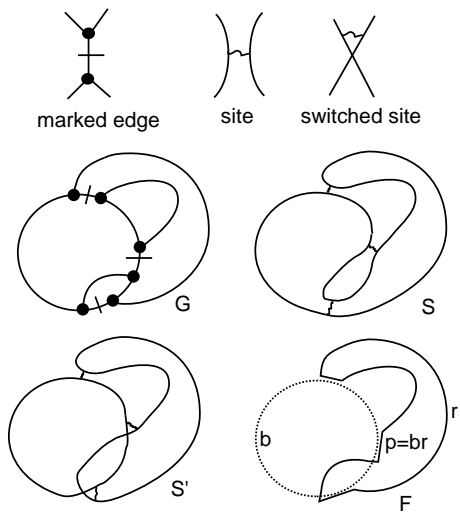}
     \end{tabular}
     \caption{\bf Graph, State, Switch and Formation.}
     \label{state}
\end{center}
\end{figure}

\begin{figure}[htb]
     \begin{center}
     \begin{tabular}{c}
     \includegraphics[width=6cm]{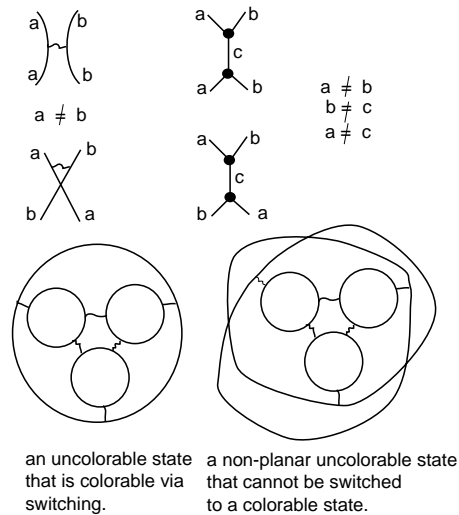}
     \end{tabular}
     \caption{\bf Crossed and Uncrossed Sites and Two Examples of States}
     \label{site}
\end{center}
\end{figure}

In this section we give a reformulation of the map theorem that is very close in spirit to the diagrammatics of low-dimensional knot theory.
Consider a collection of disjoint Jordan curves drawn in the plane. Let there be a collection of {\it sites} between distinct curves as indicated in Figures ~\ref{state} and \ref{site}.
A site is a region in the plane containing two arcs, each from one of the component curves, equipped with an indicator, marking and connecting the two curves.
A marked edge in a graph can be converted to a site by removing the edge and its nodes and connecting the local arcs as shown at the top of Figure~\ref{state}.
The indicator is to be interpreted as an instruction that the two curves at the site are to be colored differently from a set of three distinct colors $\{r,b,p\}.$ In Figure~\ref{state} we illustrate a graph $G$ 
with three edges marked that form a perfect matching (as described in the previous section). Adjacent to $G$ on the right is a state $S$ with three sites that correspond to the edges in the perfect matching.
Note that the two local arcs at a site may belong to a single component or to distinct components of the state $S.$
Below these two parts of the figure, on the left, is a diagram of a state $S'$ obtained by switching two sites of the state $S$ so that $S'$ is colorable. The extra crossings in the diagram for $S'$ show where the
switching was applied. The coloring of $S'$ is rewritten as a formation in the 
diagram on the lower right of the figure.\\ 

\begin{figure}[htb]
     \begin{center}
     \begin{tabular}{c}
     \includegraphics[width=6cm]{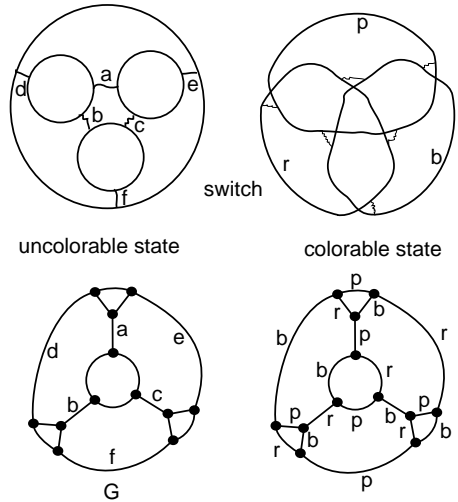}
     \end{tabular}
     \caption{\bf Coloring a Planar State and Its Graph.}
     \label{three}
\end{center}
\end{figure}

\begin{figure}[htb]
     \begin{center}
     \begin{tabular}{c}
     \includegraphics[width=6cm]{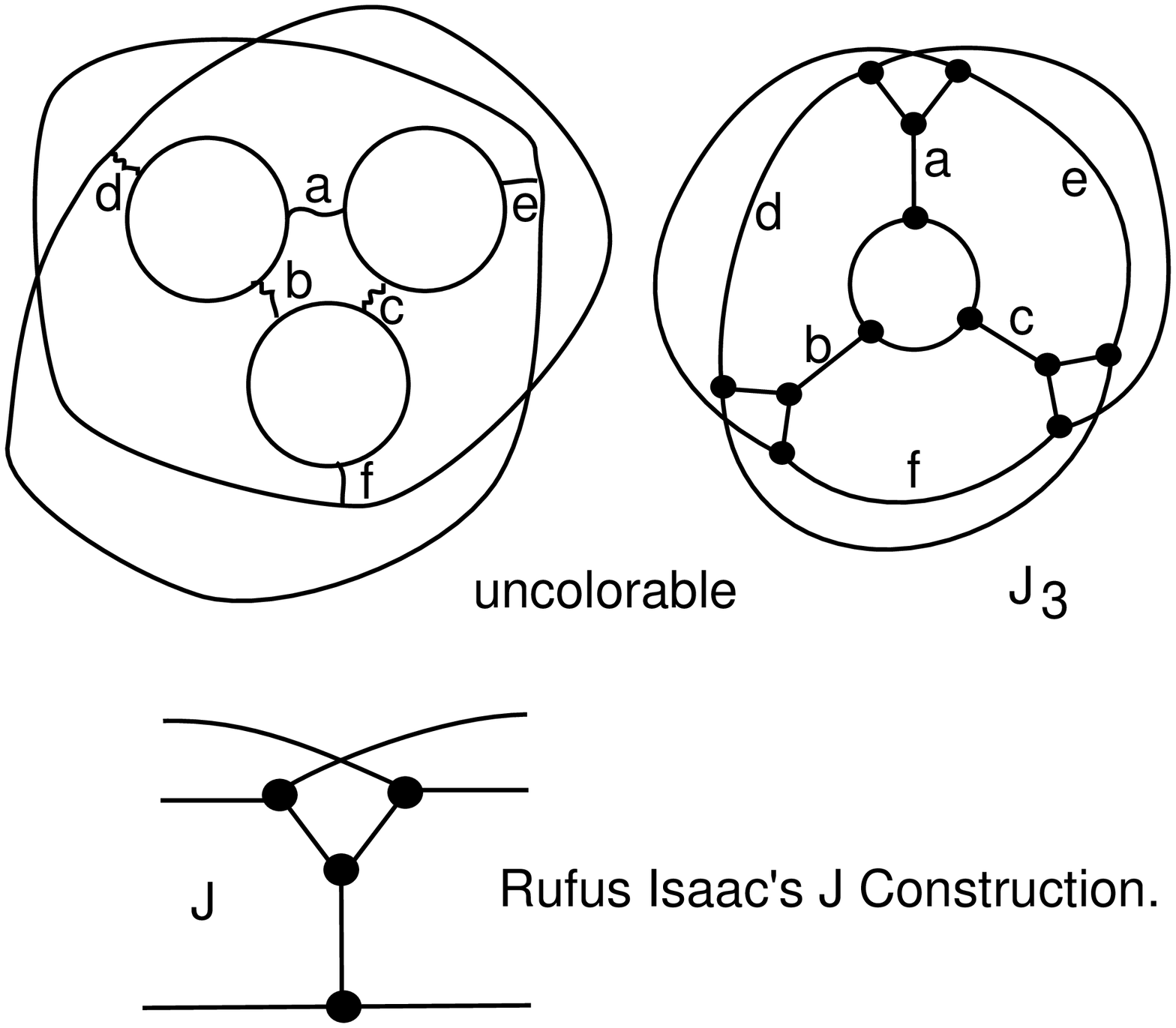}
     \end{tabular}
     \caption{\bf A NonPlanar State Corresponds to Isaac's $J_{3}.$ }
     \label{IsaacJ}
\end{center}
\end{figure}

Note that each site corresponds, as shown in Figure~\ref{state}, to an edge in a cubic graph. In fact, given a cubic graph and a perfect matching of its nodes, we can convert each pair of nodes in the matching to a site. See Figure~\ref{state} for an illustration of this remark. In the figure we have marked two edges to be so converted and, in this case, we obtain a single Jordan curve $S$ from the replacements.  A collection of Jordan curves in the plane with a chosen collection of sites will be called a {\it planar state S}. Given a planar state S, we can attempt to color its curves according to the indicators.
It may be uncolorable. For example, the state S in Figure~\ref{state} is not colorable. We introduce the local operation of {\it site switching} as show in Figure~\ref{state}. In switching a site we replace two
locally parallel arcs with two arcs that cross one another. In the resulting state we are still interested in coloring its loops. Loops that cross are allowed to have the same or different colors unless there is a state indicator forcing them to be colored differently. When we switch a site, the site indicator remains in place as shown in both Figure~\ref{state} and Figure~\ref{site}. Note that in Figure~\ref{state} we chose to switch two sites and arrive at a new state $S'$ that is colorable. This gives a coloring of the original graph $G,$ and a formation $F$ for that graph that corresponds directly to the coloring of the state $S'.$ Figure~\ref{state}, when examined should
be sufficient for the reader to understand the relationships between graphs, formations, states and sites.\\

The examples shown in Figure~\ref{site} are of interest. At the bottom left of the figure we see a state that is uncolorable but can be converted to a colorable state by switching. The state itself is not directly colorable since it consists in four mutually touching planar loops. In fact there is only one way to color this example. One must switch all of the sites as shown in Figure~\ref{three}. One then obtains
a configuration of three loops and a coloring of the corresponding graph. It is remarkable how close this example comes to being uncolorable. By some miracle it just manages to be colorable when we switch all the sites.  \\

The other example in Figure~\ref{site} is a non-planar variant of the example we just discussed. This example was drawn in an attempt convert the first example to an uncolorable and we succeeded in this attempt! In Figure~\ref{IsaacJ}  we show the graph corresponding to this state. Remarkably, the graph is $J_{3}$ the well-known uncolorable graph constructed by Rufus Isaacs \cite{Isaacs} as a circular composition
of three copies of the special graphical element $J$ shown in Figure~\ref{site}. It is not hard to see that an odd number of the $J$ elements placed in a circular configuration will always be uncolorable. These graphs are called $J_{2n+1}$ for $n = 1,2, \cdots .$ The pathway to re-discovering the Isaac's construction that we have given here is illuminating. It shows that graphs in the plane can be precariously close to uncolorability, and that a seemingly slight modification into non-planarity can produce uncolorables. Uncolorable, isthmus-free, cubic graphs were named {\em snarks} by Martin Gardner \cite{M} in association with the elusive and mysterious creature in Lewis Carroll's poem ``The Hunting of the Snark" \cite{LC}. Isaac's snarks are certainly fundamental. In fact, the graph $J_{3}$ can be 
collapsed to the minimal non-planar uncolorable, the Petersen graph. We will return to this example and discuss the Petersen graph below. \\

 We say that a planar site is an {\it isthmus} if there is a pathway from one side of the site to the other that does not cross any of the Jordan curves.  
{\it We can now give a statement of our state-calculus reformulation of the 
map theorem.}\\

\noindent {\bf Planar State Theorem.} Let $S$ be a planar state with no isthmus. Then, after switching some subset of the sites of $S$ to make a new (possibly non-planar) state $S'$, the state $S'$ is colorable with three colors.\\

\noindent {\bf Proof.} In order to prove this Theorem, we show how a planar state corresponds to a cubic planar map. This is illustrated in Figure~\ref{state} and Figure~\ref{site}. In these figures we show how  to associate a pair of cubic nodes to each 
site of a planar state. We call this operation {\it squeezing} the site. As figures show, we can squeeze a site or a switched site, obtaining the same local graphical configuration. By applying this association to each site in the state, we obtain a cubic planar map $M(S),$ that is isthmus--free exactly when the state has no isthmus as defined above.
If $S'$ is a state obtained from $S$ by switching some site, then the same map $M(S)$ results from squeezing all the sites. When the state $S'$ is colored with three colors, then the map $M(S')= M(S)$ is
colored on its edges with three colors. Conversely, if we begin with a cubic planar map $G,$ then we can associate a planar state to $G$, by choosing perfect matching of the nodes of $G$. Perfect matchings
exist for planar cubic maps, as discussed in the previous section. The perfect matching selects a special set of edges in $G$ so that every node is a unique end of one of these special nodes. By resolving the special nodes (inverse of squeezing) we obtain a planar state $S(G)$ such that $M(S(G)) = G.$ We are now in a position to see that $G$ is edge $3$--colorable if and only if $S(G)$ is switch-colorable.  
Figure~\ref{three}  illustrates how a coloring of the state $S$ induces a coloring of the map $M(S).$ Note also that we show how in squeezing a site we induce a coloring on the resulting double--Y graph
by taking the third color determined by the two distinct colors at the site. This third color is assigned to the edge created in the squeeze. Starting with a colored graph and a perfect matching we can resolve the coloring to a specific choice of switching of $S(G)$ that yields a coloring of $S(G).$ Thus we see that the problem of edge-coloring cubic maps in the plane is equivalent to the problem of coloring planar sites by 
switching. This completes the proof of the Theorem.  $\hfill\Box$
\\

\begin{figure}[htb]
     \begin{center}
     \begin{tabular}{c}
     \includegraphics[width=6cm]{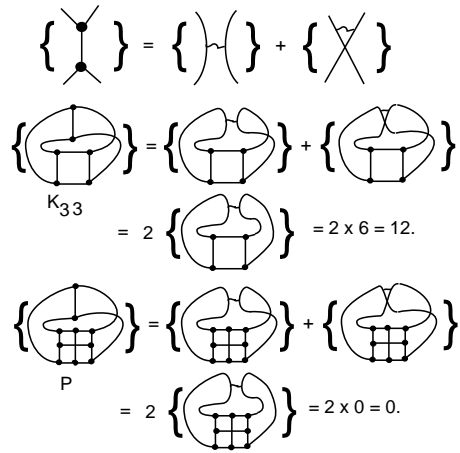}
     \end{tabular}
     \caption{\bf A Logical State Expansion, $K_{3,3}$ and Petersen Graph Examples.}
     \label{logical}
\end{center}
\end{figure}

\begin{figure}[htb]
     \begin{center}
     \begin{tabular}{c}
     \includegraphics[width=6cm]{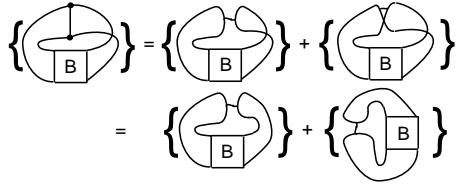}
     \end{tabular}
     \caption{\bf Expanding a General Edge.}
     \label{logicaledge}
\end{center}
\end{figure}

\begin{figure}[htb]
     \begin{center}
     \begin{tabular}{c}
     \includegraphics[width=6cm]{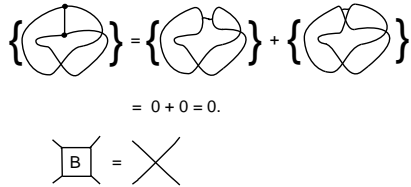}
     \end{tabular}
     \caption{\bf Dumbell from Crossover Black Box.}
     \label{dumbbox}
\end{center}
\end{figure}

Now we will use the state calculus to give a (tautological) expansion formula for colorings of arbitrary cubic maps. The formula is given below and in Figure~\ref{logical}.

$$\{ \YGlyph  \} = \{ \VGlyph \} + \{  \CGlyph  \}$$

The meaning of this equation is as follows: $\{ G \}$ denotes the number of proper three-colorings of the edges of a cubic graph $G.$ The graph $G$ is not assumed to be planar or represented in the plane.
The left-hand term of the formula represents a graph $G$ and one specific edge of $G.$ The two right hand terms of the formula represent the graph $G$ with the specific edge replaced either by a parallel connection of edges as illustrated, or by a cross-over connection of edges. In each case we have added the site-symbol, used in this section, to indicate that the two edges, in a coloring of the new graph, are to be colored with different colors. Thus the two graphs on the right-hand side of the formula have acquired a local state site of the type discussed in this section. The rest of the graphs may still have nodes. 
{\it One can define $\{ G \}$ directly as a summation over all states of the graph $G$, if we choose a perfect matching for $G$ and sum over all replacements of the edges in the perfect matching by parallel and crossed local sites.} Note that the crossed lines in the diagram for the formula are only an artifact of representation, and are not new nodes in these graphs. The truth of the formula follows from our previous discussion of the fact that a coloring of an edge will result in a double $Y$ form (as on the left-hand side of the formula) so that the two colors at the top of the $Y$ either match in parallel or match when crossed. If we are given the pattern of matched colors that are different, then we can reconstruct an edge and give it that color that is different from both of them. We shall refer to this formula as the 
{\it logical state expansion for three-coloring of cubic graphs.}\\

In Figure~\ref{logical} we give two examples of the use of this logical expansion. In the first example we expand the graph $K_{3,3}$ and find that $\{K_{3,3}\} = 12.$ Here we expand on one edge, utilize the 
symmetry of the situation to leave only one graph to consider and then see by inspection that this graph has six colorings. In the second example, we expand the Petersen graph (the reader can check that the 
initial graph is a version of the Petersen graph). In this case symmetry reduces the calculation to one graph, but it is easy to see that this graph is not colorable so that the marked edges are colored with different 
colors. Hence the Petersen graph is uncolorable.\\

In Figure~\ref{logicaledge}  we use a ``black box" $B$ to illustrate the general case of expanding one edge in a graph $G.$ Here we depict the graph with one cross-over. There is no loss of generality in this 
depiction, since the crossover could be cancelled by another one inside the black box. Note that the edges entering the black box meet the corners of the box -- these are not cubic nodes. When the graph
is expanded on one edge in this configuration we see that the two terms correspond to two ways to close the edges that emanate from the black box. If both of these closures force the same color on the marked edges, then the graph $G$ is uncolorable. In this sense, the existence of an uncolorable demands a black box that will force the transmission of a color in orthogonal directions. The black box that 
produces the Petersen graph has this property and it is very hard to imagine any black box with this property that does not produce a non-planar graph. The simplest example of an uncolorable is shown in 
Figure~\ref{dumbbox}. Here the graph is planar. It is a dumbell as in Figure~\ref{peter}, but in the form of black box and extra edge with crossing, the contents of the black box is a pair of crossed lines.
Certainly these give the simplest example of a black box that forces the same colors on the external closure arcs for both ways of making the external closure.\\

In the next section, we study the Penrose state summation for determining the number of proper colorings of a cubic graph. The Penrose formula has formal similarities with the tautological formulas of this section.\\

\section {The Penrose Formula and a Generalization to All Cubic Graphs}
\begin{figure}[htb]
     \begin{center}
     \begin{tabular}{c}
     \includegraphics[width=6cm]{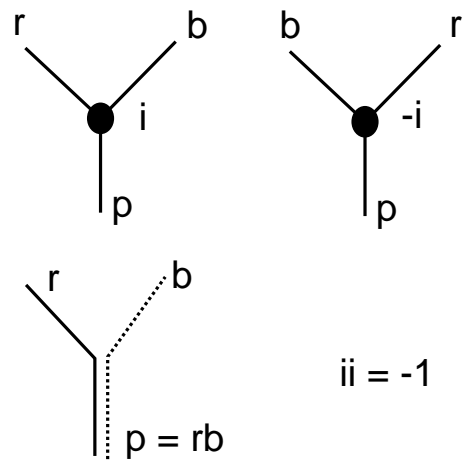}
     \end{tabular}
     \caption{\bf Node as Epsilon Tensor}
     \label{epsilon}
\end{center}
\end{figure}

Roger Penrose \cite{P} gives a formula for computing the number of proper 
edge 3-colorings of a plane cubic graph $G.$
In this formula each node is associated with the ``epsilon" tensor (for our purposes a {\it tensor} is a multi-indexed matrix. that is for each particular choice of indices the tensor returns a value in the complex numbers.)
$$P_{ijk}=\sqrt{-1}\epsilon_{ijk}$$ 
\noindent as shown in Figure~\ref{epsilon}.  Let $(xyz)$ denote the ordered list of the indices $x,y,z.$ One takes the colors from the set $\{r,b,p\}$ and the tensor $\epsilon_{ijk}$ takes value $1$
for $(ijk) = (rbp), (bpr),(prb)$ and $-1$ for $(ijk) = (rpb), (pbr), (brp).$  The tensor is $0$ when $(ijk)$ is not 
a permutation of $(rpb).$   Note that in Figure~\ref{epsilon} that $(rbp)$ corresponds to clockwise order. Since cyclic permuations of the indices do not change the value of the epsilon, the diagrammatic tensor with index labels will return $+\sqrt{-1}$ for the clockwise order and  $-\sqrt{-1}$ for anti-clockwise order. \\

One then evaluates the graph $G$ by taking the sum, over all possible color
assignments to its edges, of the products of the $P_{ijk}$ associated with its nodes. Call this 
evaluation $[G].$ \\

\noindent {\bf Remark. } This evaluation, $[G],$ is the {\it tensor contraction} of the labelling of the nodes of $G$ with the tensor $P_{ijk}.$  Any given color assignment to the edges of $G$ assigns a specific value of
$P_{ijk}$ at each node of $G$ (in this case it is $\pm \sqrt{-1}$). By definition, the tensor contraction is the summation of the products of these evaluations at the nodes where the sum is taken over all color assignments.
Note that only proper colorings of the graph contribute to this sum, since the epsilon tensor vanishes when any two colors at a node are the same. If we had labeled the nodes of the graph with other tensors then we would include their values in the product for a given coloring. Later in this section we will add such an extra tensor for a crossing induced by immersing a non-planar graph in the plane.\\

Due to properties of the epsilon tensor,  $[G]$ satisfies the identity shown below. We refer to this identity as the {\it Penrose formula}.
$$[ \YGlyph ] = [ \VDiag ] -  [ \CDiag ]$$
Along with this identity the Penrose bracket satisfies the formula below
$$[ O \, G] = 3 [G]$$
and
$$[O] = 3,$$ where $O$ denotes a Jordan curve disjoint from $G$ in the plane.\\

\begin{figure}[htb]
     \begin{center}
     \begin{tabular}{c}
     \includegraphics[width=6cm]{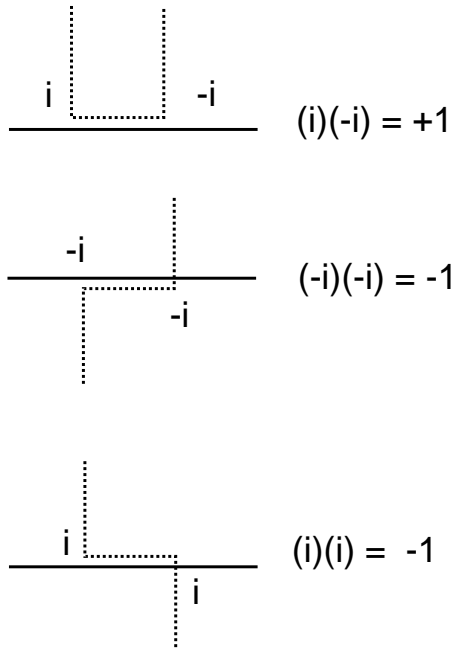}
     \end{tabular}
     \caption{\bf Bounce and Cross Contribute +1 and -1.}
     \label{bouncecross}
\end{center}
\end{figure}

\noindent {\bf Theorem (Penrose).} If $G$ is a planar cubic graph, then $[G]$, as defined above, is equal to the number of
distinct proper colorings of the edges of $G$ with three colors (so that every node sees three colors at its edges).
\bigbreak

\noindent {\bf Proof.} It follows from the above description that only proper colorings of $G$ contribute to 
the summation $[G],$ and that each such coloring contributes a product of $\pm \sqrt{-1}$ from the tensor
evaluations at the nodes of the graph. In order to see that $[G]$ is equal to the number of 
colorings for a plane graph, one must see that each such contribution is equal to $+1.$
The proof of this assertion is given in Figure~\ref{bouncecross}  with reference to Figure~\ref{epsilon} where we see that in a formation for a coloring 
each bounce contributes $+1 = -\sqrt{-1}\sqrt{-1}$ while each crossing contributes $-1.$ Since there
are an even number of crossings among the curves in the formation, it follows that the total 
product is equal to $+1$. This completes the proof of the Penrose Theorem. $\hfill\Box$
\bigbreak

\begin{figure}[htb]
     \begin{center}
     \begin{tabular}{c}
     \includegraphics[width=6cm]{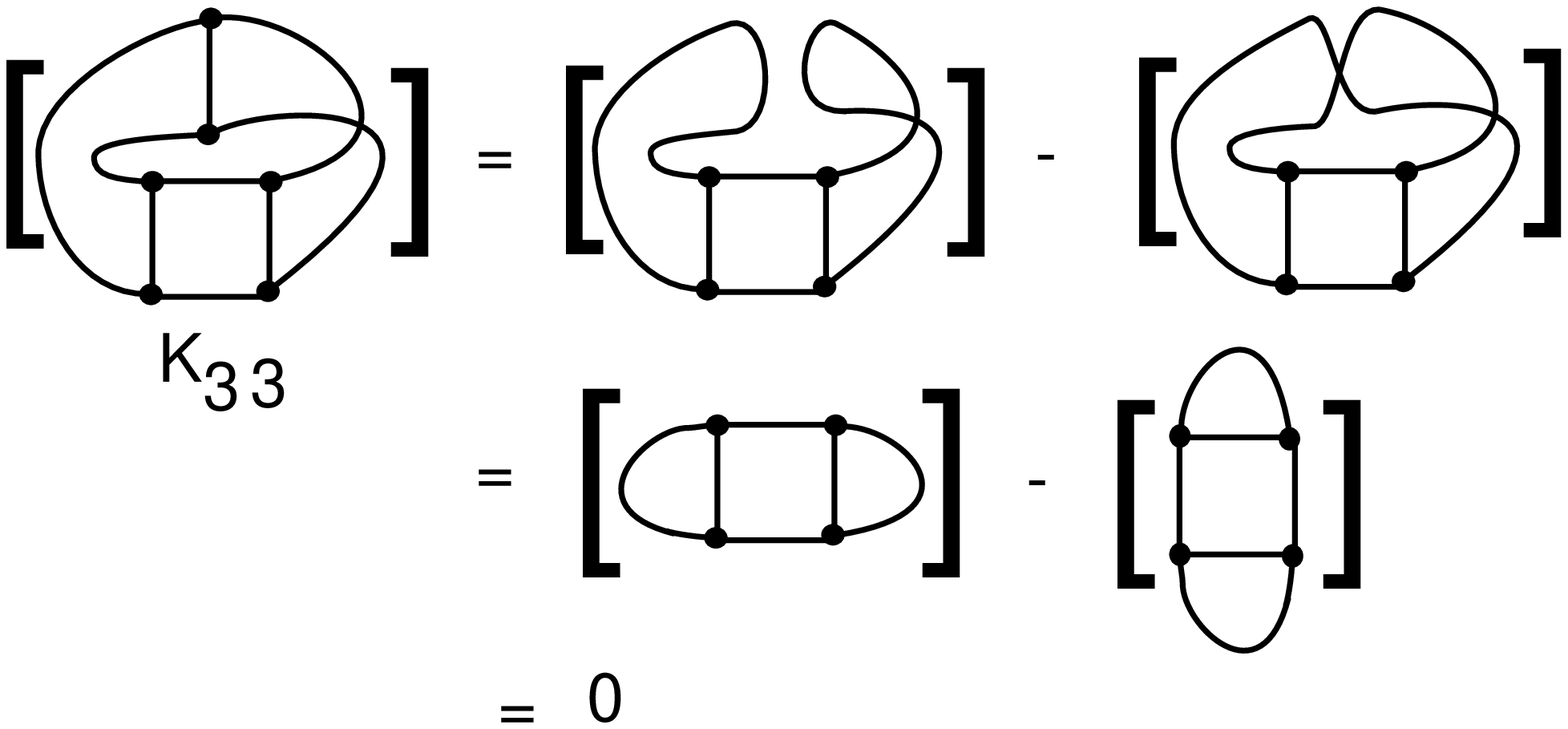}
     \end{tabular}
     \caption{\bf Penrose on K33 is Zero.}
     \label{penroseK33}
\end{center}
\end{figure}

The Penrose formula has a very interesting limitation in that, in the formulation of it that we have given so far, it only counts colors correctly for a plane graph, a graph that is given with an embedding in the plane. If we take a diagram of a non-planar graph (with crossings) and expand the Penrose bracket it will not always count the number of colorings of the graph. For example, view Figure~\ref{penroseK33} where we show that the Penrose bracket calculation of the $K_{3,3}$ graph is zero! We know from the previous section that this graph has $12$ colorings. The reason for this discrepancy should be apparent to a reader of the proof above. By using the properties of the epsilon tensor and the properties of formations as described earlier in the paper, we showed that each coloring contributed $+1$ to the state summation for a planar graph. This sign of $+1$ depended on our use of the Jordan curve theorem in counting colored Jordan curves that intersected one another an even number of times. When we use a planar diagram with crossings the number of such intersections can be odd.\\

\begin{figure}[htb]
     \begin{center}
     \begin{tabular}{c}
     \includegraphics[width=6cm]{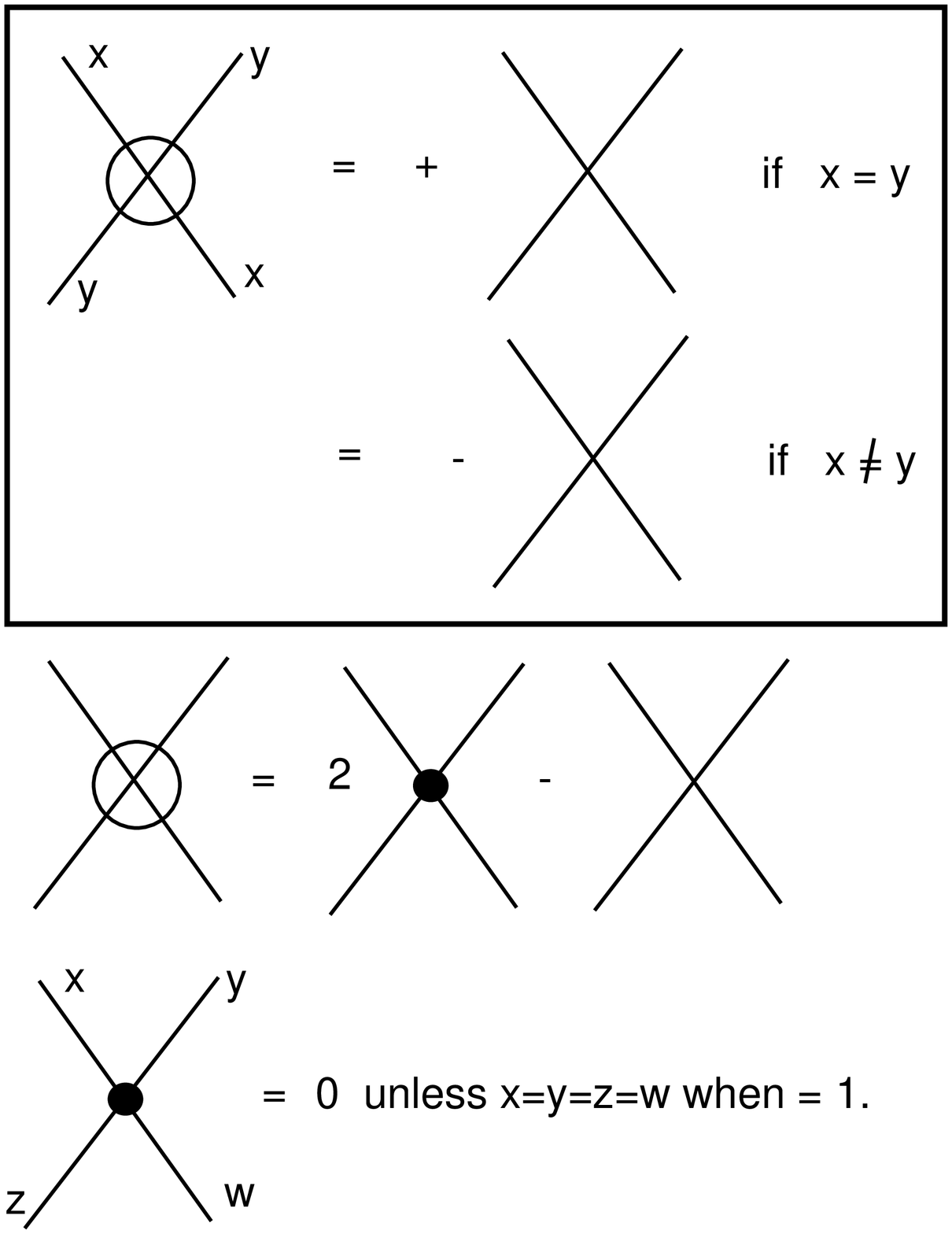}
     \end{tabular}
     \caption{\bf Crossing Tensor For Revised Penrose Bracket.}
     \label{crossingtens}
\end{center}
\end{figure}

\begin{figure}[htb]
     \begin{center}
     \begin{tabular}{c}
     \includegraphics[width=6cm]{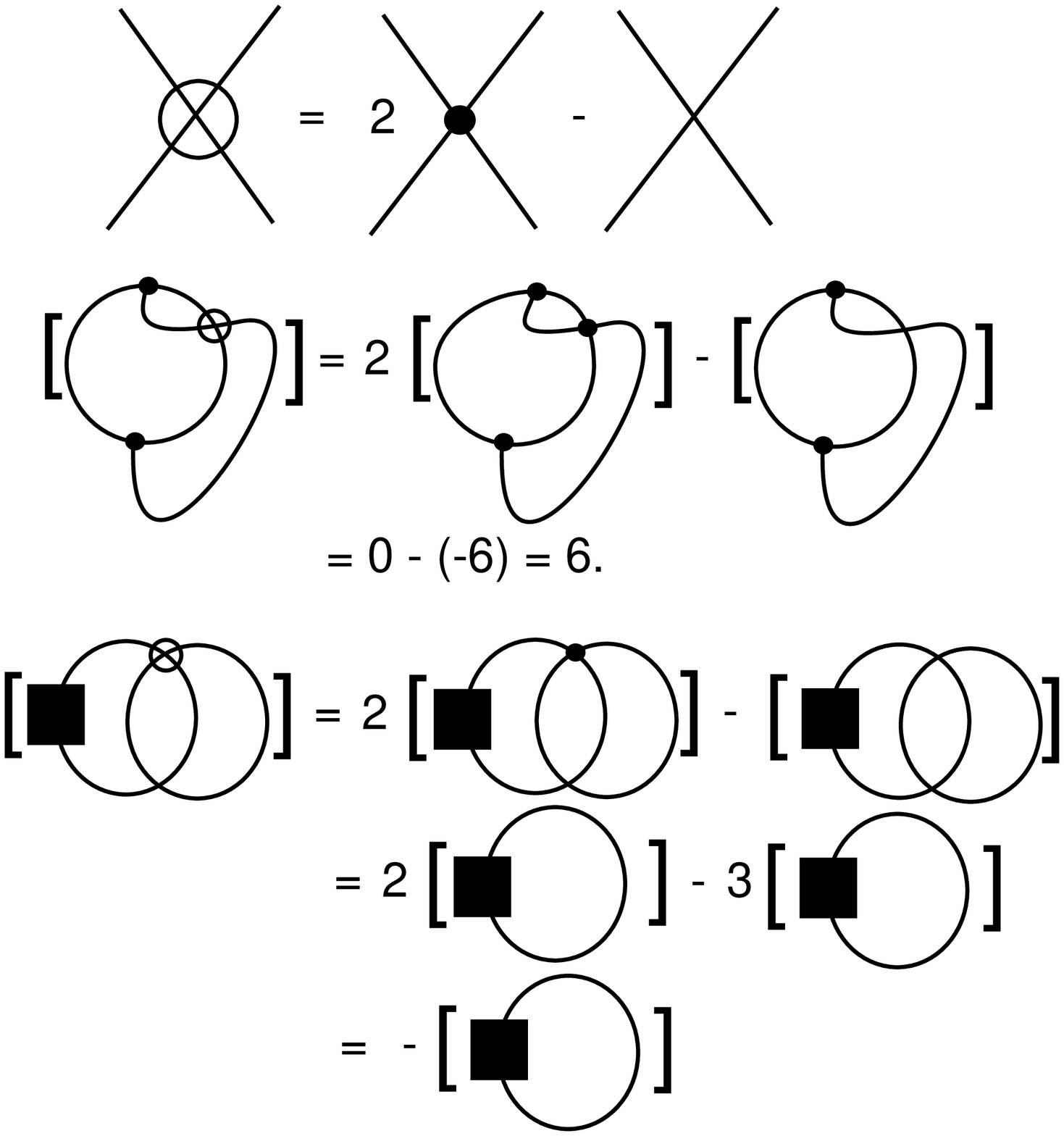}
     \end{tabular}
     \caption{\bf Crossing Tensor Formalism.}
     \label{ctensapp}
\end{center}
\end{figure}

The remarks in the previous paragraph point the way to a modification of the Penrose bracket so that it will calculate the number of colorings of any graph that is given as a diagram in the plane with crossings that are not nodes of the graph, but artifacts of the way the diagram is drawn in the plane (possibly necessary for non-planar graphs). Call these the {\em virtual crossings}.\\

Our solution to the problem is to add an extra tensor for each virtual crossing in the original diagram. The new tensor is indicated by a circle around that crossing. The circled flat crossing is a tensor depending on its four endpoints. It is zero unless the labels at the ends of a given straight segment in the crossed segments are the same. When these ends are the same we can label the crossed segments with two colors $x$ and $y,$ as in Figure~\ref{crossingtens}. Then the value of the tensor is $+1$ if $x=y$ and $-1$ if $x \neq y,$ as shown in Figure~\ref{crossingtens}. The Penrose bracket is computed just as before and it satisfies the same formulas as before. The new crossings that occur in the expansion formula are standard crossings. {\it Only the initial crossings in the diagram are circled.} It is not hard to examine the proof above and see that now every coloring will contribute $+1$ to the state sum. And thus we have the Theorem following the next remark.\\

\begin{figure}[htb]
     \begin{center}
     \begin{tabular}{c}
     \includegraphics[width=6cm]{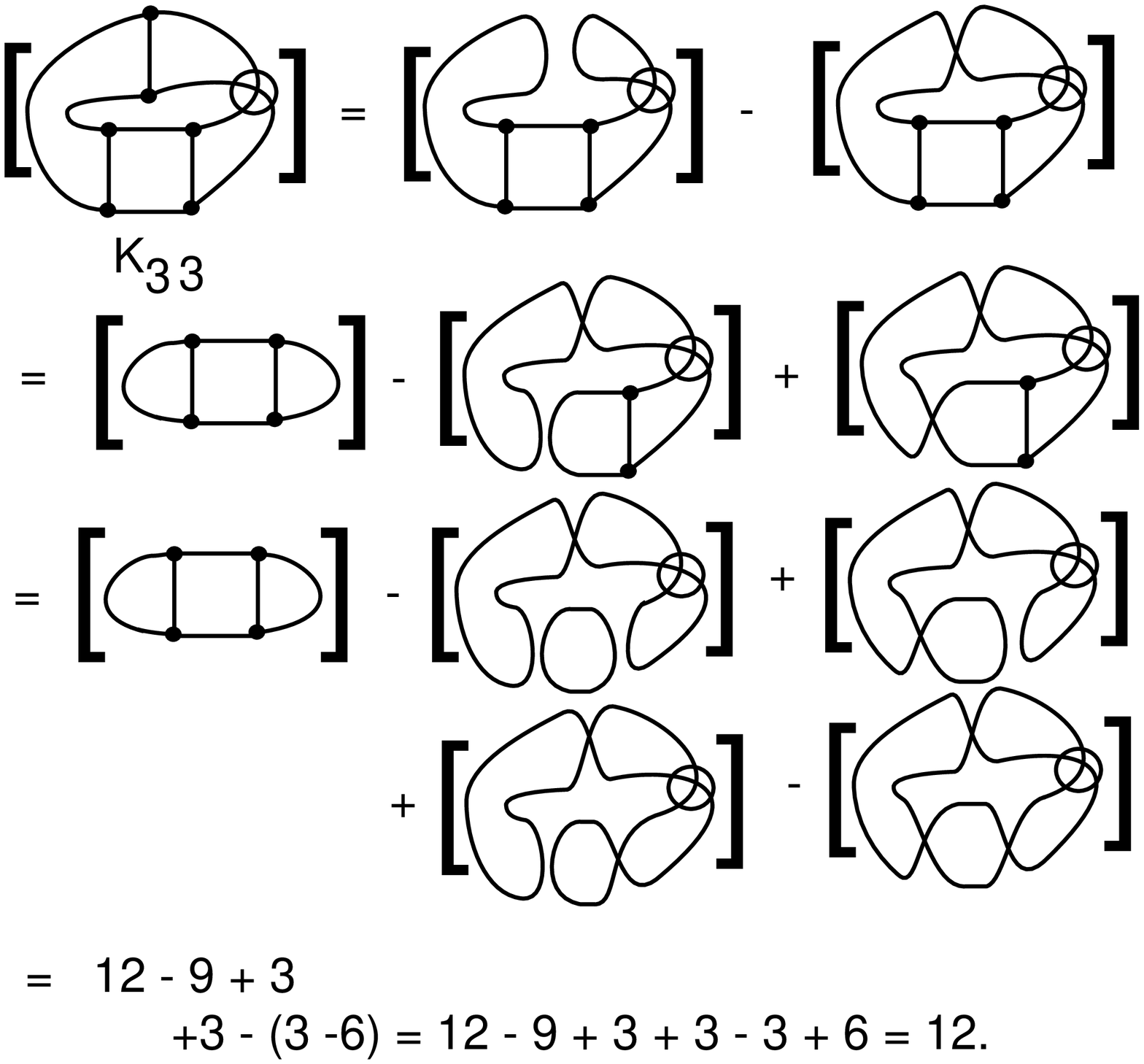}
     \end{tabular}
     \caption{\bf Revised Penrose on K33 Counts Colorings.}
     \label{revisedpenroseK33}
\end{center}
\end{figure}

\noindent {\bf Remark.} Another way to handle the circled crossing tensor is illustrated in Figures~\ref{crossingtens} and \ref{ctensapp}. Here we write formally the equation
$$[ \CCircleDiag ] = 2[ \CDotDiag ] -  [ \CDiag ].$$
The dotted crossing is regarded as a coloring node that demands that {\it all four lines incident to this node have the same color.} The ordinary crossing is taken as usual to indicate that each crossing line has a color independent of the other line. One can then expand the new bracket according to these rules and make decisions on evaluations either at the end of the process or earlier by logical considerations. 
For example, in Figure~\ref{ctensapp} we show how to evaluate a simple theta graph with a circled self-crossing and we show an example of one of many general formulas one can derive, in this case an extra circle with one circled crossing changes the sign of the evaluation. This fact could have been used in the next example that we compute directly in Figure~\ref{revisedpenroseK33}.\\

\noindent {\bf Theorem (Extending Penrose Bracket to All Cubic Graphs).} Let $G$ be any cubic graph equipped with an immersion into the plane so that the transverse crossings of interiors of edges of $G$ are circled as described above. Interpret the 
nodes of $G$ as epsilon tensors and the circled crossings as sign tensors as described above. Let $[G]$ denote the tensor contraction of the immersed graph $G$ for the three color indices. Then 
$[G]$ is equal to the number of proper three-colorings of the cubic graph $G.$ The new version of the Penrose bracket continues to satisfy the basic formulae
$$[ \YGlyph ] = [ \VDiag ] -  [ \CDiag ],$$
$$[ O \, G] = 3 [G],$$
$$[O] = 3.$$ 
The virtual crossings resulting from the immersion of $G$ are expanded by the formula
$$[ \CCircleDiag ] = 2[ \CDotDiag ] -  [ \CDiag ]$$
as described above so that the dotted crossing is regarded as a coloring node that demands that {\it all four lines incident to this node have the same color.} \\

\noindent {\bf Proof.} The proof of this result follows from the discussion preceding the statement of the Theorem. $\hfill\Box$
 \\

\begin{figure}[htb]
     \begin{center}
     \begin{tabular}{c}
     \includegraphics[width=4cm]{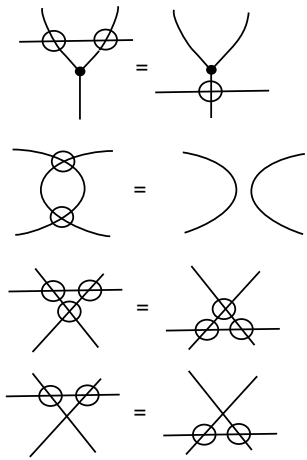}
     \end{tabular}
     \caption{\bf Topological Tensor Identities.}
     \label{toptensid}
\end{center}
\end{figure}

\begin{figure}[htb]
     \begin{center}
     \begin{tabular}{c}
     \includegraphics[width=4cm]{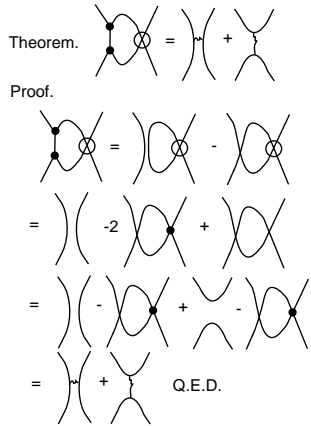}
     \end{tabular}
     \caption{\bf Basic Expansion.}
     \label{basicexpansion}
\end{center}
\end{figure}

Examine Figure~\ref{revisedpenroseK33}. In this figure we have an immersion of  $G = K_{3,3}$ with one crossing and this crossing is circled. We then expand the Penrose bracket until we have collections of circles for the second term of the first expansion. We keep the first term and rewrite it without a circled crossing because this is a self-crossing and will always have equal colors and hence 
will contribute a $+1$ in the evaluation. In the other part of the expansion to collections of circles, each final diagram can be decided, using the rule for sign of a circled crossing. Note that in the very last diagram we have two curves that cross once with a circled crossing and once with an uncircled crossing. When these two curves are colored the same the circle gives $+1$ and when they are colored unequally, the circle gives $-1.$ Thus these two circles contribute $3 -6 = -3.$ The reader will note that the initial second term adds up to zero and the total for $[G]$ is $12,$ the number of colorings of this graph.\\

\begin{figure}[htb]
     \begin{center}
     \begin{tabular}{c}
     \includegraphics[width=4cm]{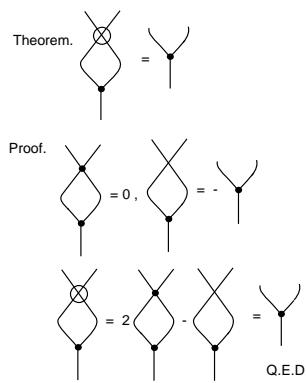}
     \end{tabular}
     \caption{\bf A Tensor Identity.}
     \label{tensid}
\end{center}
\end{figure}

Now examine Figure~\ref{tensid}. In this figure we show formally that the circled crossing placed as a twist at a cubic node does not change the evaluation of the special Penrose bracket. In this figure and
the figures to follow, we do not place brackets around the diagrams, but they are implicit in the calculations.  The point of this derivation is to show that the formal properties of the tensors that we have added
are sufficient to prove results that can also be directly deduced from the state sum definition of the extended bracket. In Figure~\ref{toptensid} we list moveablity properties for the circled crossings. These can be proven also either by direct appeal to the state summation, or by formal work with the tensors. It is useful to keep these properties in mind when working with graphical calculations.\\

Finally, in Figure~\ref{basicexpansion} we give a derivation of a basic identity for a double $Y$ form of two cubic nodes in conjunction with a circled crossing. This expands to two diagrams where we have used
our earlier notation for a state site, indicating color difference between the nearby arcs. This is the extended Penrose bracket version of a corresponding identity for the logical expansion of Section 4   shown 
in Figure~\ref{logicaledge} and using black box notation in that figure. This brings us full circle to our original logical considerations and shows how that logic is now part of the extended Penrose bracket.\\

We are happy to have a simple extension of the Penrose bracket that calculates the number of colorings of any cubic graph. Our point of view can be compared with \cite{EM} where a generalization of the Penrose bracket is made for graphs embedded in surfaces. Our state sum and abstract tensor method will be the subject of papers subsequent to the present work.\\

\end{document}